\documentclass[12pt]{amsart}
\usepackage{amssymb}
\usepackage[dvips]{graphicx}
\usepackage[dvips]{graphics}
\newtheorem{theorem}{{\bf Theorem}}
\newtheorem{lemma}{{\bf Lemma}}
\newtheorem{prop}{{\bf Proposition}}
\newtheorem{cor}{{\bf Corollary}}

\begin{document}

\title{THE VIRTUAL AND UNIVERSAL BRAIDS}

\author[Bardakov]{Valerij G. Bardakov}
\address{Sobolev Institute of Mathematics, Novosibirsk 630090, Russia}
\email{bardakov@math.nsc.ru}

%
\thanks{Authors were supported in part by the Russian Foundation
for Basic Research (grant~02--01--01118).}
\date{\today}

\begin{abstract}
We study the structure of the virtual braid group. It is shown
that the virtual braid group is a semi--direct product of the
virtual pure braid group and the symmetric group. Also, it is
shown that the virtual pure braid group is a semi--direct product
of free groups. From these results we obtain  a normal form of
words in the virtual braid group. We introduce the concept of a
universal braid group. This group contains the classical braid
group and has as its quotient groups the singular braid group,
virtual braid group, welded braid group, and classical
braid group.
\end{abstract}

\maketitle






Recently some generalizations of
classical knots and links were  defined and studied: singular links
 \cite{Vas, Bir93}, virtual links \cite{Ka, GPV} and welded links \cite{FRR}.

One of the ways to study classical links is to study the braid group.
Singular braids
 \cite{Baez, Bir93}, virtual braids \cite{Ka, Ver01}, welded
braids \cite{FRR} were defined similar to the classical braid group.
Theorem of A.~A.~ Markov \cite[Ch. 2.2]{Bir}
reduces the problem of classification of links to some algebraic
problems of the theory of braid groups.
These problems include the word problem and the conjugacy problem.
There are generalizations of  Markov theorem for singular links \cite{Ge},
virtual links, and welded links
 \cite{Kam}.

There are some different ways to solve the word problem for
the singular braid monoid and  singular braid group
 \cite{Dash, Corran, Ver}.
The solution of the word problem for the welded braid group
follows
from the fact that this group is a subgroup of the automorphism group of
the free group
\cite{FRR}. A normal form of words in the welded braid group was
constructed in
\cite{Kr}.

In this paper we study the structure of the virtual braid group
$VB_n$. Similar to the classical braid group
 $B_n$ and welded braid group $WB_n$, the group $VB_n$ contains the normal subgroup $VP_n$
which is called  {\it virtual pure braid group}. The quotient group $VB_n/VP_n$ is isomorphic
to the symmetric
group $S_n$. In the article we find generators and defining relations of $VP_n$.
Since
$VB_n$ is a semi--direct product of $VP_n$ and $S_n$, we should  study the structure of  $VP_n$.
It will be proved that $VP_n$ is representable as the following semi--direct product
$$
VP_n = V_{n-1}^* \leftthreetimes VP_{n-1} = V_{n-1}^* \leftthreetimes
(V^*_{n-2} \leftthreetimes (\ldots \leftthreetimes
(V_2^* \leftthreetimes V_1^*))\ldots),
$$
where $V_i^*$ is  some (in general infinitely generated  for $i > 1$) free subgroup of $VP_n$.
From this result it follows that there exists a normal form of words in $VB_n$.

In the last section we define the universal braid group $UB_n$ which contains the braid group
$B_n$ and has as its quotient groups the singular braid group $SG_n$,
virtual braid group $VB_n$, welded braid group $WB_n$, and braid group $B_n$.
It is known \cite{FRR} that $VB_n$ has as its quotient the group $WB_n$.
It will be proved that the quotient homomorphism maps $VP_n$ into the welded pure braid group
$WP_n$.
This homomorphism agrees with the decomposition of this
group into the semi-direct product given by Theorem~\ref{theorem2} and
by \cite{Bar, BarP}.

By  Artin theorem,  $B_n$ is embedded into the automorphism group  $\mbox{Aut}(F_n)$
of the free group $F_n$. In \cite{FRR} it was proved  that  $WB_n$ is also embedded into $\mbox{Aut}(F_n)$.
It is not known if it is true that $SG_n$ and $VB_n$ are embedded into
$\mbox{Aut}(F_n)$.

\subsection*{Acknowledgments}
I am very grateful to Joanna Kania-Bartos\-zyn\-ska, Jozef Przytycki,
Pawel Traczyk, and Bronislaw Wajnryb for organizing  of the Mini--semester on
Knot Theory (Poland, July, 2003) and for the invitation to
participate in this very interesting and very well--organized
Mini--semester.
I would also like to thank Vladimir Vershinin and Andrei Vesnin for their
interest to this work.
 Special thanks goes to the participants of the
seminar ``Evariste Galois'' at Novosibirsk State
University for their kind attention to my work.

\section{different classes of braids and their properties}

In this section we  remind (see references from the introduction) some known facts about braid
groups, singular braid monoids, virtual braid groups and welded
braid groups.

\subsection{The braid group and the group of conjugating automorphisms}

\begin{figure}[h!]
\begin{center}
\includegraphics[scale=0.30,bb= 0 0 500 500]{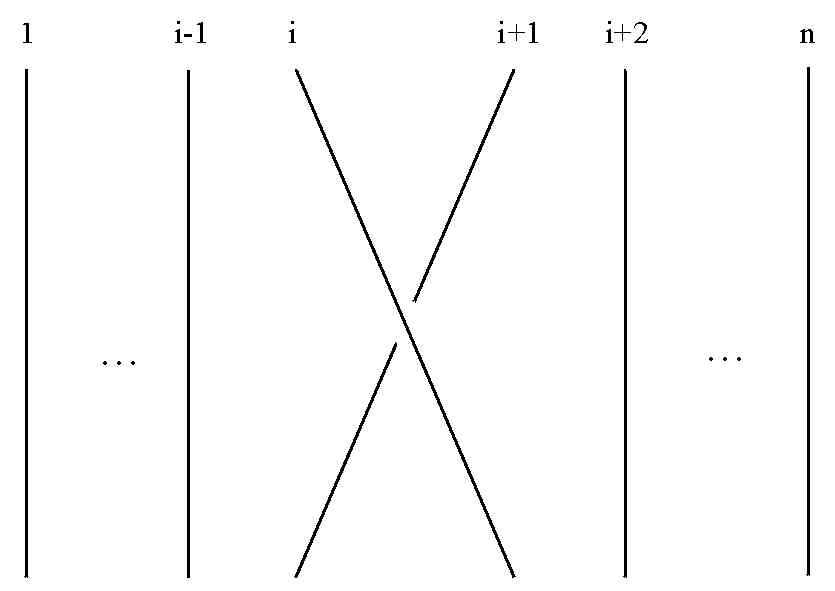}
{\caption{Geometric braids representing $\sigma_{i}$ and $\sigma_{i}^{-1}$}\label{f:1}}
 \end{center}
\end{figure}

The braid group $B_n$, $n\geq 2$, on $n$ strings can be defined as
a group generated by $\sigma_1,\sigma_2,\ldots,\sigma_{n-1}$  (see
Fig. \ref{f:1})


 \begin{figure}[h!]
 \begin{center}
 \includegraphics[scale=0.30,bb= 0 0 500 500]{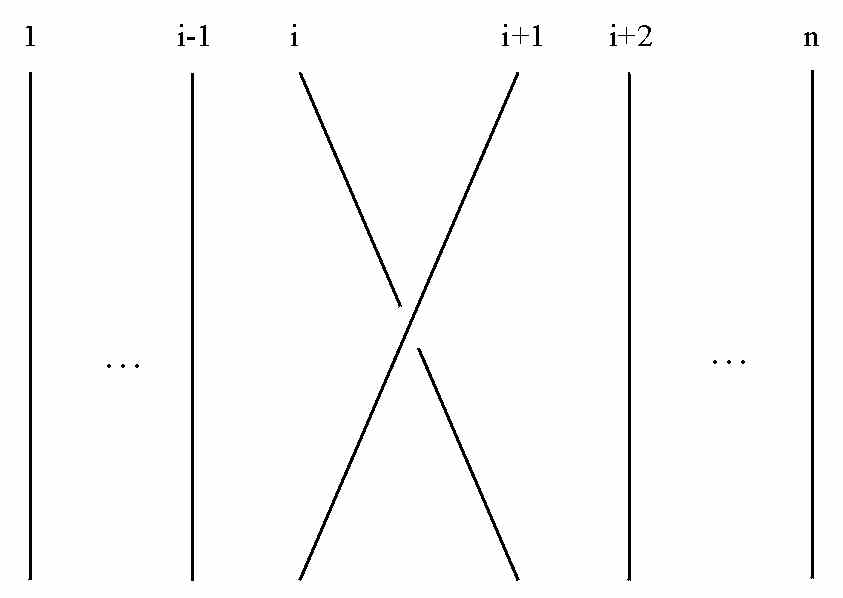}
 \end{center}
\end{figure}


with the defining relations
\begin{equation}
\sigma_i \, \sigma_{i+1} \, \sigma_i = \sigma_{i+1} \, \sigma_i \, \sigma_{i+1},~~~ i=1,2,\ldots,n-2, \label{eq1}
\end{equation}
\begin{equation}
\sigma_i \, \sigma_j = \sigma_j \, \sigma_i,~~~|i-j|\geq 2. \label{eq2}
\end{equation}

\begin{figure}[h!]
 \begin{center}
 \includegraphics[scale=0.30,bb= 0 0 500 500]{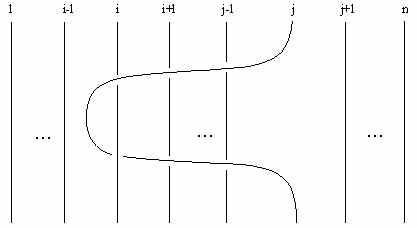}
 {\caption{The geometric braid $a_{ij}$}\label{f:3}}
 \end{center}
\end{figure}

There exists a homomorphism of $B_n$ onto the symmetric group $S_n$ on
$n$ letters. This homomorphism  maps
 $\sigma_i$ to the transposition  $(i,i+1)$, $i=1,2,\ldots,n-1$.
The kernel of this homomorphism is called
{\it pure braid group} and denoted by
$P_n$. The group $P_n$ is generated by  $a_{ij}$, $1\leq i < j\leq n$
(see Fig. \ref{f:3}).
These
generators can be expressed by the generators of
 $B_n$ as follows
$$
a_{i,i+1}=\sigma_i^2,
$$
$$
a_{ij} = \sigma_{j-1} \, \sigma_{j-2} \ldots \sigma_{i+1} \, \sigma_i^2 \, \sigma_{i+1}^{-1} \ldots
\sigma_{j-2}^{-1} \, \sigma_{j-1}^{-1},~~~i+1< j \leq n.
$$

The group $P_n$ is the semi--direct product of  the normal subgroup
$U_n$ which is a free group with free generators $a_{1n}, a_{2n},\ldots,a_{n-1,n},$
and $P_{n-1}$.
Similarly, $P_{n-1}$ is the semi--direct product of the free group
$U_{n-1}$ with free generators $a_{1,n-1}, a_{2,n-1},\ldots,a_{n-2,n-1}$ and $P_{n-2},$ and so on.
Therefore, $P_n$ is decomposable (see \cite{Mar}) into the following semi--direct product
$$
P_n=U_n\leftthreetimes (U_{n-1}\leftthreetimes (\ldots \leftthreetimes
(U_3\leftthreetimes U_2))\ldots),~~~U_i\simeq F_{i-1}, ~~~i=2,3,\ldots,n.
$$

The group $B_n$ has a faithful representation as a group of automorphisms of
 ${\rm
Aut}(F_n)$ of the free group $F_n = \langle x_1, x_2, \ldots, x_n \rangle.$
In this case the generator $\sigma_i$, $i=1,2,\ldots,n-1$, defines the automorphism
$$
\sigma_{i} : \left\{
\begin{array}{ll}
x_{i} \longmapsto x_{i} \, x_{i+1} \, x_i^{-1}, &  \\ x_{i+1} \longmapsto
x_{i}, & \\ x_{l} \longmapsto x_{l}, &  l\neq i,i+1.
\end{array} \right.
$$

By theorem of Artin \cite[Theorem 1.9]{Bir}, an automorphism $\beta $ from ${\rm
Aut}(F_n)$ lies in  $B_n$ if and only if $\beta $ satisfies to the following conditions:
$$
~~~~~1)~~ \beta(x_i) = a_i^{-1} \, x_{\pi(i)} \, a_i,~~~1\leq i\leq n,
$$
$$
2)~~ \beta(x_1x_2 \ldots x_n)=x_1x_2 \ldots x_n,
$$
where $\pi $ is a permutation from $S_n$ and $a_i\in F_n$.

An automorphism is called a  {\it conjugating automorphism} (or a permutation--conjugating
automorphism according to the terminology from
 \cite{FRR}) if it satisfies to condition 1).
The group of conjugating automorphisms $C_n$ is generated by $\sigma_i$
and automorphisms $\alpha_i$, $i=1, 2, \ldots,$ $n-1$, where
$$
\alpha_{i} : \left\{
\begin{array}{ll}
x_{i} \longmapsto x_{i+1}, & \\
x_{i+1} \longmapsto x_{i},  &\\
x_{l} \longmapsto x_{l}, &  l\neq i, i+1.
\end{array} \right.
$$
It is not hard to see that the automorphisms $\alpha_i$ generate
the symmetric group $S_n$ and, hence, satisfy the following relations
\begin{equation}
\alpha_i \, \alpha_{i+1} \, \alpha_i = \alpha_{i+1} \, \alpha_i \, \alpha_{i+1},~~~ i=1,2,\ldots,n-2, \label{eq3}
\end{equation}
\begin{equation}
\alpha_i \, \alpha_j = \alpha_j \, \alpha_i,~~~ |i-j|\geq 2, \label{eq4}
\end{equation}
\begin{equation}
\alpha_i^2 = 1,~~~ i=1,2,\ldots,n-1. \label{eq5}
\end{equation}
The group $C_n$ is defined by relations (\ref{eq1})--(\ref{eq2}) of $B_n$,
relations (\ref{eq3})--(\ref{eq5})
of $S_n$, and the mixed relations (see \cite{FRR, Sav})
\begin{equation}
\alpha_i \, \sigma_j = \sigma_j \, \alpha_i,~~~ |i-j|\geq 2, \label{eq6}
\end{equation}
\begin{equation}
\sigma_i  \, \alpha_{i+1} \, \alpha_i = \alpha_{i+1}  \, \alpha_i  \, \sigma_{i+1},~~~ i=1,2,\ldots,n-2, \label{eq7}
\end{equation}
\begin{equation}
\sigma_{i+1} \, \sigma_{i} \, \alpha_{i+1} = \alpha_{i} \, \sigma_{i+1} \, \sigma_{i},~~~ i=1,2,\ldots,n-2. \label{eq8}
\end{equation}

If we consider the group  generated by automorphisms $\varepsilon_{ij}$,
$1\leq i\neq j \leq n$, where
$$
\varepsilon_{ij} : \left\{
\begin{array}{ll}
x_{i} \longmapsto x_{j}^{-1} \, x_i \, x_j, &  i\neq j, \\
x_{l} \longmapsto x_{l}, &  l\neq i,
\end{array} \right.
$$
then we get the group of {\it basis--conjugating automorphisms} $Cb_n$.
The elements of $Cb_n$ satisfy  condition 1) for the identical permutation $\pi $, i.~e.,
map each generator $x_i$ to the conjugating element.
J.~McCool
 \cite{Mac} proved that $Cb_n$ is defined by the relations
(from here different letters stand for different indices)
\begin{equation}
\varepsilon_{ij} \, \varepsilon_{kl} = \varepsilon_{kl} \, \varepsilon_{ij},
 \label{eq9}
\end{equation}
\begin{equation}
\varepsilon_{ij} \, \varepsilon_{kj} = \varepsilon_{kj} \, \varepsilon_{ij},
 \label{eq10}
\end{equation}
\begin{equation}
(\varepsilon_{ij} \, \varepsilon_{kj}) \, \varepsilon_{ik} = \varepsilon_{ik} \,
(\varepsilon_{ij} \,
\varepsilon_{kj}).
\label{eq11}
\end{equation}

The group $C_n$ is representable as the semi--direct product: $C_n = Cb_n \leftthreetimes S_n$,
where $S_n$ is generated by the automorphisms $\alpha_1, \alpha_2, \ldots, \alpha_{n-1}$.
The following equalities are true
(see \cite{Sav}):
$$
\varepsilon_{i,i+1} = \alpha_i \, \sigma^{-1}_i,~~~
\varepsilon_{i+1,i} = \sigma^{-1}_i \, \alpha_i,
$$
$$
\varepsilon_{ij} = \alpha_{j-1} \, \alpha_{j-2} \ldots \alpha_{i+1} \, \varepsilon_{i,
i+1} \,
\alpha_{i+1} \ldots \alpha_{j-2} \, \alpha_{j-1}~~~i <j,
$$
$$
\varepsilon_{ji} = \alpha_{j-1} \, \alpha_{j-2} \ldots \alpha_{i+1} \,
\alpha_i \,
\varepsilon_{i, i+1} \, \alpha_i \,
\alpha_{i+1} \ldots \alpha_{j-2} \, \alpha_{j-1}~~~i <j.
$$

The structure of $Cb_n$ was studied in \cite{Bar, BarP}. There it was proved that
  $Cb_n$, $n\geq 2$,
is decomposable into the semi--direct product
$$
Cb_n = D_{n-1}\leftthreetimes (D_{n-2}\leftthreetimes (\ldots
\leftthreetimes (D_2\leftthreetimes D_1))\ldots ),
$$
of subgroups $D_i,$ $i=1,2,\ldots,n-1,$ generated by
 $\varepsilon_{i+1,1},$ $\varepsilon_{i+1,2},$
$\ldots,\varepsilon_{i+1,i},$
$\varepsilon_{1,i+1},$ $\varepsilon_{2,i+1},$ $\ldots,\varepsilon_{i,i+1}$.
The elements $\varepsilon_{i+1,1},$ $\varepsilon_{i+1,2},$ $
\ldots,\varepsilon_{i+1,i}$  generate a free group of rank $i$,
elements
$\varepsilon_{1,i+1},$ $\varepsilon_{2,i+1},$ $\ldots,\varepsilon_{i,i+1}$
generate a free abelian group of rank $i$.

The  pure braid group $P_n$ is contained in  $Cb_n$ and the generators of
$P_n$ can be written in the form
$$
a_{i,i+1}=\varepsilon_{i,i+1}^{-1} \, \varepsilon_{i+1,i}^{-1},~~~~~i=1,2,\ldots,n-1,
$$
$$
a_{ij} = \varepsilon_{j-1,i} \, \varepsilon_{j-2,i} \ldots
\varepsilon_{i+1,i} \,
(\varepsilon_{ij}^{-1} \, \varepsilon_{ji}^{-1}) \, \varepsilon_{i+1,i}^{-1} \ldots
\varepsilon_{j-2,i}^{-1} \, \varepsilon_{j-1,i}^{-1} =
$$
$$
= \varepsilon_{j-1,j}^{-1} \, \varepsilon_{j-2,j}^{-1} \ldots
\varepsilon_{i+1,j}^{-1} \,
(\varepsilon_{ij}^{-1} \, \varepsilon_{ji}^{-1}) \, \varepsilon_{i+1,j} \ldots
\varepsilon_{j-2,j} \, \varepsilon_{j-1,j},~~~~~2 \leq i+1 < j \leq n.
$$

\vskip 20pt

\subsection{The singular braid monoid.}

{\it The Baez--Birman monoid}
\cite{Baez, Bir93} or {\it the singular braid monoid} $SB_n$ is generated
(as monoid) by elements $\sigma_i,$ $\sigma_i^{-1}$, $\tau_i$, $i = 1, 2, \ldots, n-1$.
The elements $\sigma_i,$ $\sigma_i^{-1}$  generate the braid group
$B_n$. The generators $\tau_i$  satisfy the defining relations
\begin{equation}
\tau_i \, \tau_j = \tau_j \, \tau_i, ~~~|i - j| \geq 1, \label{eq12}
\end{equation}
other relations are mixed:
\begin{equation}
\tau_{i} \, \sigma_{j} = \sigma_{j} \, \tau_{i}, ~~~|i - j| \geq 1, \label{eq13}
\end{equation}
\begin{equation}
\tau_{i}  \, \sigma_{i} = \sigma_{i} \, \tau_{i},~~~ i=1,2,\ldots,n-1,  \label{eq14}
\end{equation}
\begin{equation}
\sigma_{i} \, \sigma_{i+1} \, \tau_i = \tau_{i+1} \, \sigma_{i}  \, \sigma_{i+1},~~~ i=1,2,\ldots,n-2,
 \label{eq15}
 \end{equation}
 \begin{equation}
\sigma_{i+1}  \, \sigma_{i} \, \tau_{i+1} = \tau_{i} \,
\sigma_{i+1} \, \sigma_{i}, ~~~ i=1,2,\ldots,n-2.
 \label{eq16}
\end{equation}

In the work \cite{FKR} it was proved that the singular braid monoid $SB_n$ is embedded into
the group $SG_n$ which is called the  {\it
singular braid group} and has the same defining relations as $SB_n$.

\subsection{The virtual braid group and welded braid group}

\begin{figure}[h!]
 \begin{center}
 \includegraphics[scale=0.30,bb= 0 0 500 500]{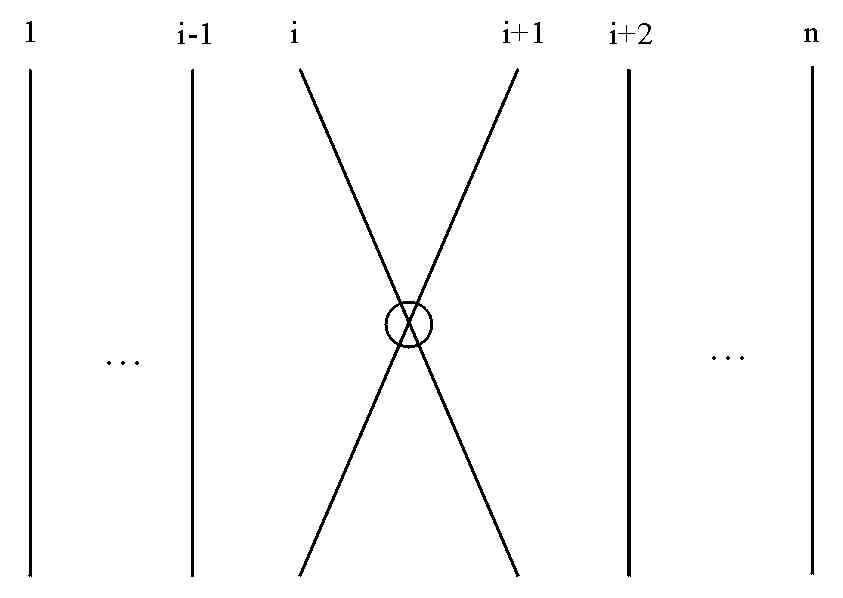}
 {\caption{The geometric virtual braid $\rho_{i}$}\label{f:2}}
 \end{center}
 \end{figure}

The virtual braid group $VB_n$ was introduced in \cite{Ka}.
In \cite{Ver01} it was found more short system of defining relations (see below).
The group $VB_n$ is generated by $\sigma_i$, $\rho_i$, $i = 1, 2, \ldots, n-1$ (see Fig. \ref{f:2}).

The elements $\sigma_i$ generate the braid group $B_n$ with defining relations
(\ref{eq1})--(\ref{eq2}) and the elements $\rho_i$ generate the symmetric group
$S_n$ which is defined by the relations
\begin{equation}
\rho_{i}  \, \rho_{i+1} \, \rho_{i} = \rho_{i+1} \, \rho_{i} \, \rho_{i+1},~~~ i=1,2,\ldots,n-2,
\label{eq17}
\end{equation}
\begin{equation}
\rho_{i} \, \rho_{j} = \rho_{j} \, \rho_{i},~~~|i-j| \geq 1,
\label{eq18}
\end{equation}
\begin{equation}
\rho_{i}^2 = 1~~~ i=1,2,\ldots,n-1.
\label{eq19}
\end{equation}
Other relations are mixed:
\begin{equation}
\sigma_{i} \, \rho_{j} = \rho_{j} \, \sigma_{i},~~~|i-j| \geq 1,
\label{eq20}
\end{equation}
\begin{equation}
\rho_{i} \, \rho_{i+1}  \, \sigma_{i} = \sigma_{i+1} \, \rho_{i}  \, \rho_{i+1},~~~
i=1,2,\ldots,n-2.
\label{eq21}
\end{equation}

Note that the last relation is equivalent to the following relation:
$$
\rho_{i+1} \, \rho_{i}  \, \sigma_{i+1} = \sigma_{i} \, \rho_{i+1}  \, \rho_{i}.
$$

In the work \cite{GPV} it was proved that the  relations
$$
\rho_{i} \, \sigma_{i+1}  \, \sigma_{i} = \sigma_{i+1} \, \sigma_{i}
\, \rho_{i+1},~~~~~
\rho_{i+1} \, \sigma_{i}  \, \sigma_{i+1} = \rho_{i} \, \sigma_{i+1}
\, \sigma_{i}.
$$
are not fulfilled in $VB_n$.

In the work \cite{FRR} it was introduced the welded braid group $WB_n$.
This group is generated by $\sigma_i$, $\alpha_i$, $i=1, 2, \ldots, n-1$. The elements
$\sigma_i$ generate the braid group $B_n$. The elements
$\alpha_i$  generate the symmetric group $S_n$ and the following mixed relations hold
\begin{equation}
\alpha_{i}  \, \sigma_{j} = \sigma_{j} \, \alpha_{i}, ~~~|i - j| \geq 1, \label{eq22}
\end{equation}
\begin{equation}
\sigma_{i} \, \alpha_{i+1} \, \alpha_i = \alpha_{i+1} \, \alpha_{i} \, \sigma_{i+1},
~~~ i=1,2,\ldots,n-2, \label{eq23}
 \end{equation}
 \begin{equation}
\sigma_{i+1} \, \sigma_{i} \, \alpha_{i+1} = \alpha_{i} \, \sigma_{i+1} \, \sigma_{i},~~~ i=1,2,\ldots,n-2.
 \label{eq24}
\end{equation}

In the work \cite{FRR} it was proved that $WB_n$ is isomorphic to the group of conjugating automorphisms
 $C_n$.

Comparing the defining relations of $VB_n$ with the defining relations of $WB_n$,  we see that
$WB_n$ can be obtained from $VB_n$ by adding some new relation.
Therefore, there exists a homomorphism
$$
\varphi_{VW} : VB_n \longrightarrow WB_n,
$$
taking  $\sigma_i$ to $\sigma_i$ and $\rho_i$ to $\alpha_i$ for all $i$.
Hence, $WB_n$
is the homomorphic image of $VB_n$.

In \cite{FRR} it was proved  that the following relation (symmetric to (\ref{eq23}))
$$
\sigma_{i+1} \, \alpha_{i} \, \alpha_{i+1} = \alpha_{i} \, \alpha_{i+1} \, \sigma_{i},
$$
is true in $WB_n$. But the following relation is not fulfilled
$$
\alpha_{i+1}  \, \sigma_{i} \, \sigma_{i+1} = \sigma_{i} \, \sigma_{i+1}
\, \alpha_{i}.
$$

In \cite{Ver01} it was constructed the linear representations of $VB_n$
and $WB_n$ by matrices from $\mbox{GL}_n(\mathbb{Z}[t, t^{-1}])$ which
continue  the well known Burau representation.
 The linear representation of $C_n \simeq
WB_n$ it was constructed in \cite{BarP}. This representation
continue (with some conditions on parameters) the known Lawrence--Krammer
representation.

\section{Generators and defining relations of the virtual pure braid group}

In this section we  introduce a virtual pure braid group and find its generators and defining relations.

Define the map
$$
\nu : VB_n \longrightarrow S_n
$$
of $VB_n$ onto the symmetric group $S_n$ by actions on generators
$$
\nu(\sigma_i) = \nu(\rho_i) = \rho_i, ~~~ i = 1, 2, \ldots, n-1,
$$
where $S_n$ is the group generated by $\rho_i$.
The kernel $\mbox{ker}(\nu)$ of this map is called the
{\it virtual pure braid group} and denoted by $VP_n$.
It is clear that $VP_n$ is a normal subgroup of index $n!$ of $VB_n$.
Moreover, since  $VP_n \bigcap S_n = e$ and $VB_n = VP_n \cdot S_n$,
then
$VB_n = VP_n \leftthreetimes S_n$, i.~e., the virtual pure braid group is the semi--direct
product of $VP_n$ and $S_n$.

\begin{figure}[h!]
 \begin{center}
 \includegraphics[scale=0.30,bb= 0 0 500 500]{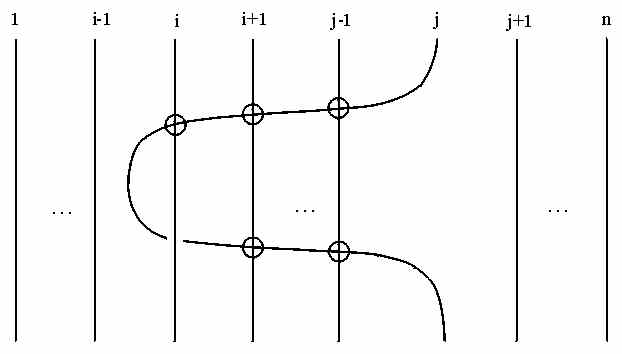}
 {\caption{The geometric virtual braid $\lambda_{ij}$ $(1 \leq i < j \leq n)$}\label{f:5}}
 \end{center}
 \end{figure}

\begin{figure}[h!]
 \begin{center}
 \includegraphics[scale=0.30,bb= 0 0 500 500]{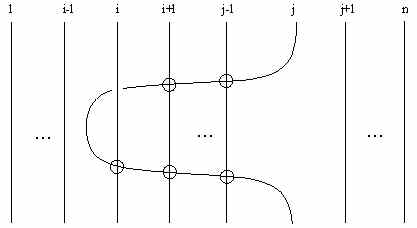}
{\caption{The geometric virtual braid $\lambda_{ji}$ $(1 \leq i < j \leq n)$}\label{f:6}}
 \end{center}
 \end{figure}

Define the following elements
$$
\lambda_{i,i+1} = \rho_i \, \sigma_i^{-1},~~~
\lambda_{i+1,i} = \rho_i \, \lambda_{i,i+1} \, \rho_i = \sigma_i^{-1} \, \rho_i,
~~~i=1, 2, \ldots, n-1,
$$
$$
\lambda_{ij} = \rho_{j-1} \, \rho_{j-2} \ldots \rho_{i+1} \, \lambda_{i,i+1} \, \rho_{i+1}
\ldots \rho_{j-2} \, \rho_{j-1},
$$
$$
\lambda_{ji} = \rho_{j-1} \, \rho_{j-2} \ldots \rho_{i+1} \, \lambda_{i+1,i} \, \rho_{i+1}
\ldots \rho_{j-2} \, \rho_{j-1}, ~~~1 \leq i < j-1 \leq n-1.
$$
Obviously,  all these elements belong to $VP_n$
and have the following geometric interpretation
(Fig. \ref{f:5}, \ref{f:6})

The next lemma hold

\begin{lemma} \label{lemma1}
Let  $1 \leq i < j \leq n$.  The following conjugating rules are fulfilled in $VB_n$:\\
\smallskip
1) for $k < i-1$ or $i < k < j-1$ or $k > j$
$$
\rho_k  \, \lambda_{ij} \, \rho_k = \lambda_{ij},~~~~~\rho_k \, \lambda_{ji} \, \rho_k =
\lambda_{ji};
$$
\\
2) $
\rho_{i-1} \,  \lambda_{ij} \,  \rho_{i-1} = \lambda_{i-1,j},
~~~\rho_{i-1}  \, \lambda_{ji} \,  \rho_{i-1} =
\lambda_{j,i-1};\\
$
\smallskip
3) for $i < j-1$
$$
\begin{array}{ll}
\rho_i  \, \lambda_{i,i+1} \,  \rho_i = \lambda_{i+1,i}, & \rho_i \,  \lambda_{ij} \,  \rho_i =
\lambda_{i+1,j}, \\
\rho_i \,  \lambda_{i+1,i}  \, \rho_i = \lambda_{i,i+1}, & \rho_i  \, \lambda_{ji} \,  \rho_i =
\lambda_{j,i+1}; \\
\end{array}
$$
\smallskip
4) for $i+1 < j$
$$
\rho_{j-1} \,  \lambda_{ij}  \, \rho_{j-1} = \lambda_{i,j-1},
~~~~~\rho_{j-1} \,  \lambda_{ji}  \, \rho_{j-1} =
\lambda_{j-1,i};
$$
\smallskip
5)
$
\rho_{j}  \, \lambda_{ij}  \, \rho_{j} = \lambda_{i,j+1},~~~\rho_{j} \,  \lambda_{ji} \,  \rho_{j} =
\lambda_{j+1,i}.
$
\end{lemma}

\begin{proof}
We consider only the rules containing $\lambda_{ij}$ for $i < j$
(the remaining rules can be considered analogously). Recall that
$$
\lambda_{ij} = \rho_{j-1}  \, \rho_{j-2} \ldots \rho_{i+1}  \, \lambda_{i,i+1}  \,  \rho_{i+1}
\ldots \rho_{j-2}  \, \rho_{j-1}.
$$
If $k < i - 1$ or $k > j$ then $\rho_k$ is permutable with
$\rho_i, \rho_{i+1}, \ldots, \rho_{j-1}$
in view of relation (\ref{eq18})  and with $\sigma_i$ in view of relation (\ref{eq20}).
Hence, $\rho_k$ is permutable with $\lambda_{ij}$.

Let $i < k < j-1$. Then
$$
\rho_k  \, \lambda_{ij}  \, \rho_{k} = \rho_{k}  \, (\rho_{j-1} \ldots \rho_{k+2} \,  \rho_{k+1}
 \, \rho_{k} \ldots
\rho_{i+1}  \, \lambda_{i,i+1} \,  \rho_{i+1} \ldots \rho_{k} \,  \rho_{k+1}  \, \rho_{k+2}
\ldots \rho_{j-1})  \, \rho_k.
$$
Permuting $\rho_k$ to $\lambda_{i,i+1}$ while it is possible, we get
$$
\rho_{j-1}  \, \ldots \rho_{k+2} \,  (\rho_{k}  \, \rho_{k+1}  \, \rho_{k}) \ldots
\rho_{i+1}  \, \lambda_{i,i+1} \,  \rho_{i+1} \ldots (\rho_{k} \,  \rho_{k+1}  \, \rho_{k})
 \,  \rho_{k+2}
\ldots \rho_{j-1}.
$$
Using the relation $\rho_{k}  \, \rho_{k+1} \,  \rho_{k} =
\rho_{k+1} \,  \rho_{k} \,  \rho_{k+1}$,
rewrite the last formula as follows:
$$
\rho_{j-1} \,  \ldots \rho_{k+2} \,  \rho_{k+1}  \, \rho_{k} \,  (\rho_{k+1} \,  \rho_{k-1} \ldots
\rho_{i+1} \,  \lambda_{i,i+1} \,  \rho_{i+1} \ldots \rho_{k-1} \,  \rho_{k+1}) \times
$$
$$
\times \rho_{k}
 \, \rho_{k+1} \,  \rho_{k+2}
\ldots \rho_{j-1}
= \rho_{j-1} \ldots \rho_{k}  \, (\rho_{k+1}  \, \lambda_{i,k} \,  \rho_{k+1})   \, \rho_{k}
 \ldots \rho_{j-1}.
$$
In view of the case considered earlier, we have
$$
\rho_{k+1}  \, \lambda_{ik} \,  \rho_{k+1} = \lambda_{ik}
$$
and, hence,
$$
\rho_{j-1} \ldots \rho_{k}  \, (\rho_{k+1}  \, \lambda_{ik} \,  \rho_{k+1}) \,  \rho_{k}
 \ldots \rho_{j-1} =
\lambda_{ij}.
$$
Thus, the first rule from 1) is proven.

2) Consider
$$
\rho_{i-1} \,  \lambda_{ij} \,  \rho_{i-1} = \rho_{i-1} \,  (\rho_{j-1} \,  \rho_{j-2} \ldots
\rho_{i+1} \,  \lambda_{i,i+1}  \, \rho_{i+1} \ldots \rho_{j-2}
 \, \rho_{j-1})  \, \rho_{i-1}.
$$
Using relation (\ref{eq18}), let as permute $\rho_{i-1}$ to $\lambda_{i,i+1}$
as long as it is possible.
We get
\begin{equation}
\rho_{i-1}  \, \lambda_{ij}  \, \rho_{i-1} = \rho_{j-1} \ldots \rho_{i+2} \,  \rho_{i+1}
(\rho_{i-1} \,  \lambda_{i,i+1} \,  \rho_{i-1}) \,  \rho_{i+1}
\rho_{i+2} \ldots \rho_{j-2}.
 \label{eq27}
\end{equation}
The expression in the brackets can be rewritten in the following form
$$
\rho_{i-1}  \, \lambda_{i,i+1} \,  \rho_{i-1} = \rho_{i-1}  \, \rho_{i}  \, \sigma_i^{-1}
 \,  \rho_{i-1}
= \rho_{i-1}  \, \rho_{i}  \, \sigma_i^{-1} \,  \rho_{i-1}
\rho_{i}  \, \rho_{i}.
$$
Using the relation
$\sigma_i^{-1} \,  \rho_{i-1} \,  \rho_{i} = \rho_{i-1} \,  \rho_i \,  \sigma_{i-1}^{-1}$
(it follows from (\ref{eq21})) and (\ref{eq18}), (\ref{eq19}),
we obtain
$$
\rho_{i-1} \,  \rho_{i} \,  (\sigma_{i}^{-1} \,  \rho_{i-1} \,  \rho_{i}) \,  \rho_{i} =
 \rho_{i-1}  \, ( \rho_i \, \rho_{i-1}  \, \rho_i) \,  \sigma_{i-1}^{-1} \,  \rho_{i} =
$$
$$
 = (\rho_{i-1} \,  \rho_{i-1})
  \,  \rho_i  \, \rho_{i-1} \,
\sigma_{i-1}^{-1}  \, \rho_{i} = \rho_{i}  \, \lambda_{i-1,i} \,  \rho_{i}.
$$
Then from (\ref{eq27}) we obtain
$$
\rho_{i-1}  \, \lambda_{ij} \,  \rho_{i-1} = \lambda_{i-1,j}.
$$
Therefore, the desired relations are proven.

3) The first formula follows from the definitions of $ \lambda_{i,i+1}$ and $ \lambda_{i+1,i}$.
Let us consider
$$
\rho_{i} \,  \lambda_{ij} \,  \rho_{i} = \rho_{i}  \, (\rho_{j-1}  \, \rho_{j-2} \ldots
\rho_{i+1} \,  \lambda_{i,i+1} \,  \rho_{i+1} \ldots \rho_{j-2} \,
\rho_{j-1})  \, \rho_{i}.
$$
Permuting  $\rho_i$ to $\lambda_{i,i+1}$ while it is possible, we obtain
$$
\rho_{i}  \, \lambda_{ij}  \, \rho_{i} = \rho_{j-1} \ldots \rho_{i+2} \,  (\rho_{i} \,
\rho_{i+1} \,  \lambda_{i,i+1} \,  \rho_{i+1}  \, \rho_{i}) \,
\rho_{i+2}  \, \ldots \rho_{j-1}.
$$
Rewrite the expression in the brackets  as follows
$$
\rho_{i} \,  \rho_{i+1} \,  \lambda_{i,i+1} \,  \rho_{i+1} \,  \rho_i = \rho_{i}
 \,  \rho_{i+1} \,  \rho_{i} \,
(\sigma_i^{-1} \,  \rho_{i+1}  \, \rho_{i}) = \rho_{i}  \, \rho_{i+1} \,  (\rho_{i} \,
\rho_{i+1} \,  \rho_{i}) \,  \sigma_{i+1}^{-1} =
$$
$$
= \rho_{i}  \, \rho_{i+1}  \, \rho_{i+1} \,
\rho_{i} \,  \rho_{i+1} \,  \sigma_{i+1}^{-1} = \rho_{i+1}  \, \sigma_{i+1}^{-1}.
$$
Hence,
$$
\rho_{i}  \, \lambda_{ij} \,  \rho_{i} = \rho_{j-1} \ldots \rho_{i+2} \,
(\rho_{i+1} \,
\sigma_{i+1}^{-1}) \,  \rho_{i+2} \ldots \rho_{j-1} =
\lambda_{i+1,j}.
$$
Therefore, the desired relations are proven.

4) follows from the relation $\rho_{j-1}^2 = e$ and the definition of $\lambda_{ij}$.

5) is an immediate consequence of the definition of $\lambda_{ij}$.
\end{proof}

\begin{cor}\label{cor1}
The group $S_n$ acts by conjugation  on the set
$\{ \lambda_{kl} ~ \vert 1 ~ \leq k \neq l \leq n \}.$ This action is transitive.
\end{cor}

In view of Lemma \ref{lemma1}, the subgroup $\langle \lambda_{kl} ~ \vert ~ 1 \leq k \neq l \leq n \rangle$
of $VP_n$ is normal in $VB_n$. Let us prove that this group coincides with  $VP_n$ and
let us find its generators and defining relations. For this purpose we  use
the Reidemeister--Schreier method (see, for example, \cite[Ch. 2.2]{KMS}).

Let $m_{kl} = \rho_{k-1}  \, \rho_{k-2} \ldots \rho_l$ for $l < k$ and $m_{kl} = 1$
in other cases. Then the set
$$
\Lambda_n = \left\{ \prod\limits_{k=2}^n m_{k,j_k} \vert 1 \leq j_k
\leq k \right\}
$$
is a Schreier set of coset representatives of $VP_n$ in $VB_n$.

\begin{theorem} \label{theorem1}
The group $VP_n$ admits a presentation with the generators $\lambda_{kl},$ $1 \leq k \neq l \leq
n$,
and the defining relations:
\begin{equation}
\lambda_{ij} \,  \lambda_{kl} = \lambda_{kl}  \, \lambda_{ij},
 \label{eq29}
\end{equation}
\begin{equation}
\lambda_{ki} \,  (\lambda_{kj} \,  \lambda_{ij}) = (\lambda_{ij} \,  \lambda_{kj}) \,  \lambda_{ki},
 \label{eq28}
\end{equation}
where  distinct letters stand for distinct indices.
\end{theorem}

\begin{proof}
Define the map $^- : VB_n \longrightarrow \Lambda_n$ which takes an element
$w \in VB_n$
into the representative $\overline{w}$ from $\Lambda_n$. In this case the element
$w \overline{w}^{-1}$ belong to $VP_n$. By Theorem 2.7 from  \cite{KMS}
the group $VP_n$ is generated by
$$
s_{\lambda, a} = \lambda a \cdot (\overline{\lambda a})^{-1},
$$
where $\lambda$ run the set $\Lambda_n$ and $a$ run the set of generators of
$VB_n$.

It is easy to establish that $s_{\lambda, \rho_i} = e$ for all representatives $\lambda $
and generators
$\rho_i$. Consider the generators
$$
s_{\lambda, \sigma_i} = \lambda \sigma_i \cdot (\overline{\lambda
\rho_i})^{-1}.
$$
For $\lambda =e$ we get $s_{e,\sigma_i} = \sigma_i \rho_i = \lambda_{i,i+1}^{-1}$.
Note that $\lambda \rho_i$ is equal to  $\overline{\lambda \rho_i}$
in $S_n$. Therefore,
$$
s_{\lambda, \sigma_i} = \lambda (\sigma_i \rho_i) \lambda^{-1}.
$$
From Lemma \ref{lemma1} it follows  that each generator $s_{\lambda, \sigma_i}$
is equal to some
$\lambda_{kl}$, $1 \leq k \neq l \leq n$.
By Corollary  \ref{cor1}, the inverse statement is also true,
i.~e.,
each element $\lambda_{kl}$ is equal to some generator $s_{\lambda, \sigma_i}$. The first part
of the theorem is proven.

To find  defining relations of $VP_n$ we define
a rewriting process $\tau $. It allows us to rewrite a word which is written in the generators
of $VB_n$ and presents an element in $VP_n$ as a word in the generators of $VP_n$.
Let us associate to the reduced word
$$
u = a_1^{\varepsilon_1} \, a_2^{\varepsilon_2} \ldots
a_{\nu}^{\varepsilon_{\nu}},~~~\varepsilon_l = \pm 1,~~~a_l \in
\{\sigma_1, \sigma_2, \ldots, \sigma_{n-1}, \rho_1, \rho_2, \ldots, \rho_{n-1}
\},
$$
the word
$$
\tau(u) = s_{k_1,a_1}^{\varepsilon_1} \,  s_{k_2,a_2}^{\varepsilon_2}
\ldots s_{k_{\nu},a_{\nu}}^{\varepsilon_{\nu}}
$$
in the generators of $VP_n$, where $k_j$ is a representative of the ($j-1$)th
initial segment
of the word $u$ if $\varepsilon_j = 1$ and $k_j$ is a representative of the $j$th
initial segment of
the word $u$ if
$\varepsilon_j = -1$.

By \cite[Theorem 2.9]{KMS}, the group $VP_n$ is defined by relations
$$
r_{\mu,\lambda} = \tau (\lambda  \, r_{\mu} \,  \lambda^{-1}),~~~\lambda \in
\Lambda_n,
$$
where $r_{\mu}$ is the defining relation of $VB_n$.

Denote by
$$
r_1 = \sigma_i  \, \sigma_{i+1}  \, \sigma_i \,  \sigma_{i+1}^{-1} \,  \sigma_i^{-1}
 \, \sigma_{i+1}^{-1}
$$
the first relation of $VB_n$. Then
$$
r_{1,e} = \tau(r_1) = s_{e,\sigma_i} \,  s_{\overline{\sigma_i},\sigma_{i+1}} \,
s_{\overline{\sigma_i \sigma_{i+1}} ,\sigma_{i}} \,
s_{\overline{\sigma_i \sigma_{i+1} \sigma_i \sigma_{i+1}^{-1}} ,\sigma_{i+1}}^{-1} \,
s_{\overline{\sigma_i \sigma_{i+1} \sigma_i \sigma_{i+1}^{-1} \sigma_i^{-1}} ,\sigma_{i}}^{-1} \,
s_{\overline{r_1} ,\sigma_{i+1}}^{-1} =
$$
$$
= \lambda_{i,i+1}^{-1}  \, (\rho_i  \,  \lambda_{i+1,i+2}^{-1} \,  \rho_i) \,
(\rho_i  \,  \rho_{i+1} \,  \lambda_{i,i+1}^{-1} \,  \rho_{i+1} \,  \rho_i) \times
$$
$$
\times (\rho_{i+1}  \, \rho_{i}  \, \lambda_{i+1,i+2}  \, \rho_{i}  \, \rho_{i+1}) \,
(\rho_{i+1}  \, \lambda_{i,i+1} \,  \rho_{i+1}) \,
\lambda_{i+1,i+2}.
$$
Using the conjugating rules from Lemma \ref{lemma1},
we get
$$
r_{1,e} = \lambda_{i,i+1}^{-1}  \, \lambda_{i,i+2}^{-1} \,  \lambda_{i+1,i+2}^{-1} \,
\lambda_{i,i+1}  \, \lambda_{i,i+2} \,  \lambda_{i+1,i+2}.
$$
Therefore, the following relation
$$
\lambda_{i,i+1}  \, (\lambda_{i,i+2} \,  \lambda_{i+1,i+2}) = (\lambda_{i+1,i+2} \,
 \lambda_{i,i+2}) \,
 \lambda_{i,i+1}
$$
is fulfilled in $VP_n$.
The Remaining relations $r_{1,\lambda}$, $\lambda \in \Lambda_n$, can be obtained from this
relation using conjugation by $\lambda^{-1}$.
Bu the formulas from Lemma
\ref{lemma1}, we obtain  relations (\ref{eq28}).

Let us consider the next relation of $VB_n$:
$$
r_2 = \sigma_i  \, \sigma_j  \, \sigma_i^{-1}  \, \sigma_j^{-1},~~~|i - j| > 1.
$$
For it we have
$$
r_{2,e} = \tau(r_2) = s_{e,\sigma_i}  \, s_{\overline{\sigma_i},\sigma_{j}} \,
s_{\overline{\sigma_i \sigma_{j} \sigma_i^{-1}} ,\sigma_{i}}^{-1} \,
s_{\overline{r_2} ,\sigma_{j}}^{-1} =
$$
$$
= \lambda_{i,i+1}^{-1} \, \lambda_{j,j+1}^{-1} \, \lambda_{i,i+1} \,
\lambda_{j,j+1}.
$$
Hence, the relation
$$
\lambda_{i,i+1}  \, \lambda_{j,j+1} = \lambda_{j,j+1} \, \lambda_{i,i+1},
~~~|i - j| > 1
$$
holds in $VP_n$.
Conjugating this relation by all representatives from $\Lambda_n$, we obtain  relations
 (\ref{eq29}).

Let us prove that only trivial relations follow from all other relations of $VB_n$.
It is evident for
relations (\ref{eq17})--(\ref{eq19}) defining the group $S_n$ because
$s_{\lambda,\rho_i} = e$ for all
$\lambda \in \Lambda_n$ and $\rho_i$.

Consider the mixed relation (\ref{eq21}) (relation (\ref{eq20})
can be considered  similarly):
$$
r_3 = \sigma_{i+1} \, \rho_i \, \rho_{i+1} \, \sigma_i^{-1} \, \rho_{i+1} \,
\rho_i.
$$
Using the rewriting process, we get
$$
r_{3,e} = \tau(r_3) = s_{e,\sigma_{i+1}} \,
s_{\overline{\sigma_{i+1} \rho_i \rho_{i+1} \sigma_i^{-1}}
,\sigma_{i}}^{-1} =
$$
$$
= \lambda_{i+1,i+2}^{-1} \, (\rho_i  \, \rho_{i+1} \, \lambda_{i,i+1} \, \rho_{i+1} \, \rho_i) =
 \lambda_{i+1,i+2}^{-1} \, \lambda_{i+1,i+2} = e.
$$
Therefore, $VP_n$ is defined by relations (\ref{eq29})
--(\ref{eq28}).
\end{proof}

\section{The structure of the virtual braid group}

From the definition of $VP_n$ and Lemma \ref{lemma1} it follows that
$VB_n = VP_n \leftthreetimes S_n$, i.~e., $VB_n$ is the splittable extension
of the group $VP_n$ by  $S_n$. Consequently, we have to study the structure
of the virtual pure braid group $VP_n$. Let us define the subgroups
$$
V_i = \langle \lambda_{1,i+1}, \lambda_{2,i+1}, \ldots, \lambda_{i,i+1};
\lambda_{i+1,1}, \lambda_{i+1,2}, \ldots, \lambda_{i+1,i} \rangle ,~~~i=1,
2, \ldots, n-1,
$$
of $VP_n$. Each $V_i$ is a subgroup of $VP_{i+1}$. Let $V_i^*$ be
the normal closure of $V_i$ in $VP_{i+1}$. The following theorem is
the main result of this section.

\begin{theorem} \label{theorem2}
The group $VP_n$, $n \geq 2$, is representable as the semi--direct product
$$
VP_n = V_{n-1}^* \leftthreetimes VP_{n-1} = V_{n-1}^* \leftthreetimes
(V^*_{n-2} \leftthreetimes (\ldots \leftthreetimes
(V_2^* \leftthreetimes V_1^*))\ldots),
$$
where $V_1^*$ is a free group of rank $2$ and $V_i^*$, $i=2,
3, \ldots, n-1,$ are free infinitely generated subgroups.
\end{theorem}

Let us prove the theorem by induction on $n$. For $n = 2$, we have
$$
VP_2 = V_1 = V_1^*
$$
and, by Theorem \ref{theorem1}, the group $V_1$ is free generated by
 $\lambda_{12}$ and $\lambda_{21}$.

To make the general case more clear  consider the case $n=3$.

\subsection{The structure of $VP_3$.} By Theorem \ref{theorem1},
the group $VP_3$ is generated by subgroups $V_1$, $V_2$ and defined by the relations
$$
\lambda_{12} \, (\lambda_{13} \,  \lambda_{23}) = (\lambda_{23} \,
\lambda_{13})  \, \lambda_{12},~~~~~
\lambda_{21}  \, (\lambda_{23} \,  \lambda_{13}) = (\lambda_{13} \,
\lambda_{23})  \, \lambda_{21},
$$
$$
\lambda_{13}  \, (\lambda_{12} \,  \lambda_{32}) = (\lambda_{32} \,
\lambda_{12})  \, \lambda_{13},~~~~~
\lambda_{31} \,  (\lambda_{32} \,  \lambda_{12}) = (\lambda_{12} \,
\lambda_{32})  \, \lambda_{31},
$$
$$
\lambda_{23}  \, (\lambda_{21} \,  \lambda_{31}) = (\lambda_{31} \,
\lambda_{21})  \, \lambda_{23},~~~~~
\lambda_{32} \,  (\lambda_{31}  \, \lambda_{21}) = (\lambda_{21} \,
\lambda_{31})  \, \lambda_{32}.
$$

From these relations we obtain  the next lemma.

\begin{lemma} \label{lemma2}
In $VP_3$ the following equalities hold: \\
1)
$$
\begin{array}{lll}
\lambda_{13}^{\lambda_{12}} = \lambda_{32}^{\lambda_{12}} \,  \lambda_{13} \,
\lambda_{32}^{-1},&
\lambda_{31}^{\lambda_{12}} = \lambda_{32} \,  \lambda_{31} \,
\lambda_{32}^{-\lambda_{12}},&
\lambda_{23}^{\lambda_{12}} = \lambda_{13}  \, \lambda_{23} \,
\lambda_{32} \,  \lambda_{13}^{-1} \,  \lambda_{32}^{-\lambda_{12}},\\
 & & \\
\lambda_{13}^{\lambda_{12}^{-1}} = \lambda_{32}^{-1}  \, \lambda_{13} \,
\lambda_{32}^{\lambda_{12}^{-1}},&
\lambda_{31}^{\lambda_{12}^{-1}} =
\lambda_{32}^{-\lambda_{12}^{-1}} \,
\lambda_{31} \lambda_{32},&
\lambda_{23}^{\lambda_{12}^{-1}} =
\lambda_{32}^{-\lambda_{12}^{-1}} \,
\lambda_{13}^{-1}  \, \lambda_{32} \,  \lambda_{23} \,  \lambda_{13},\
\end{array}
$$
2)
$$
\begin{array}{lll}
\lambda_{23}^{\lambda_{21}} = \lambda_{31}^{\lambda_{21}} \,  \lambda_{23} \,
\lambda_{31}^{-1},&
\lambda_{32}^{\lambda_{21}} = \lambda_{31}  \, \lambda_{32} \,
\lambda_{31}^{-\lambda_{21}},&
\lambda_{13}^{\lambda_{21}} = \lambda_{23}  \, \lambda_{13} \,
\lambda_{31} \,  \lambda_{23}^{-1} \,  \lambda_{31}^{-\lambda_{21}},\\
 & & \\
\lambda_{23}^{\lambda_{21}^{-1}} = \lambda_{31}^{-1} \,  \lambda_{23} \,
\lambda_{31}^{\lambda_{21}^{-1}},&
\lambda_{32}^{\lambda_{21}^{-1}} =
\lambda_{31}^{-\lambda_{21}^{-1}} \,
\lambda_{32} \,  \lambda_{31},&
\lambda_{13}^{\lambda_{21}^{-1}} =
\lambda_{31}^{-\lambda_{21}^{-1}} \,
\lambda_{23}^{-1}  \, \lambda_{31} \,  \lambda_{13} \,  \lambda_{23},\\
\end{array}
$$
where $a^b$ stand for $b^{-1} a b$.
\end{lemma}

\begin{proof}
The first and  second relations from 1) immediately follow from the
third and  forth relations of $VP_3$ (see the relations before the lemma).
Similarly, the first and second relations from 2)
immediately follow from the fifth and  sixth relations of
$VP_3$.

Further, from the first and second relations of $VP_3$ we
obtain
$$
\lambda_{23}^{\lambda_{1,2}} = \lambda_{13}  \, \lambda_{23} \,
\lambda_{13}^{-\lambda_{12}},~~~
\lambda_{13}^{\lambda_{21}} = \lambda_{23}  \, \lambda_{13} \,
\lambda_{23}^{-\lambda_{21}}.
$$
Using the proved formulas for  $\lambda_{13}^{\lambda_{12}}$ and
$\lambda_{23}^{\lambda_{21}}$, we get the third formulas from 1) and 2) respectively.

The formulas for conjugation by $\lambda_{12}^{-1}$ and $\lambda_{21}^{-1}$ can be
obtained analogously.
\end{proof}

Note that there exists an epimorphism
$$
\varphi_3 : VP_3 \longrightarrow VP_2,
$$
which takes the generators of $V_2 = \langle \lambda_{13}, \lambda_{23}, \lambda_{31},
\lambda_{32} \rangle $
into the unit and fixes the generators of $V_1 = \langle \lambda_{12}, \lambda_{21} \rangle $.
The kernel of this epimorphism is the normal closure of
$V_2$ in $VP_3$, i.~e., $\mbox{ker}(\varphi_3) = V_2^*$.

Let $u$ be the empty word or a reduced word  beginning with
non-zero power of $\lambda_{12}$ and representing an element from $V_1$.
Let $\lambda_{32}(u) = \lambda_{32}^u$ $= u^{-1} \, \lambda_{32} \, u$. We call this element
{\it the reduced power
of the generator} $\lambda_{32}$ with the power $u$. Analogously, if $v$
is the empty word or a reduced word  beginning with
non-zero power of $\lambda_{21}$ and representing an element from
$V_1$, then we put $\lambda_{31}(v) = \lambda_{13}^v$ and  call this element
{\it the reduced power
of generator} $\lambda_{31}$ with the power $v$.

\begin{lemma} \label{lemma3}
The group $V_2^*$ is a free group with generators $\lambda_{13}$, $\lambda_{23}$
and all reduced powers of  $\lambda_{31}$ and $\lambda_{32}$.
\end{lemma}

\begin{proof}
To prove the lemma we can use the Reidemeister--Shreier method, but it is
simpler to
use the definitions of normal closure and semi-direct product.
Evidently, the group
$V_2^*$ is generated by the elements
$$
\lambda_{13}^w,~~~ \lambda_{23}^w,~~~ \lambda_{31}^w,~~~
\lambda_{32}^w,~~~w \in V_1.
$$
In view of Lemma \ref{lemma2},  it is sufficient to take from these elements
only $\lambda_{13}$, $\lambda_{23}$ and all reduced powers of the generators $\lambda_{31}$
and $\lambda_{32}$.

The freedom of $V_2^*$ follows from the representation of $VP_3$ as the semi--direct product.
Indeed, since  $V_1 \bigcap V_2^* = e$, $V_1 V_2^* = VP_3$, then
$VP_3 = V_2^* \leftthreetimes V_1$.
In this case the defining relations of $VP_3$ are equivalent to the conjugating rules
from Lemma \ref{lemma2}. Therefore, all relations define the action
of the group $V_1$ on the group $V_2^*$. Since there are not other relations,  this means
that  $V_1$ and $V_2^*$ are free groups.
\end{proof}

As a consequence of this Lemma, we obtain the normal form of words in $VP_3$.
Any element $w$ from $VP_3$ can be written in the form
$w = w_1 w_2$,
where $w_1$ is a reduced word over the alphabet $\{ \lambda_{12}^{\pm 1}, \lambda_{21}^{\pm 1} \}$
and $w_2$ is a reduced word over the alphabet
$\{ \lambda_{13}^{\pm 1}, \lambda_{23}^{\pm 1}, \lambda_{31}(u)^{\pm 1},
\lambda_{32}(v)^{\pm 1} \}$, where $\lambda_{31}(u)$, $\lambda_{32}(v)$ are reduced powers
of the generators $\lambda_{31}$ and $\lambda_{32}$ respectively.

\subsection{The proof of Theorem \ref{theorem2}}
Now, we introduce the following notation.
By $\lambda_{ij}^*$ denote any
$\lambda_{ij}$ or $\lambda_{ji}$ from $VP_n$.

\begin{lemma} \label{lemma4}
For every $n \geq 2$ there exists a homomorphism
$$
\varphi : VP_n \longrightarrow VP_{n-1}
$$
which takes the generators $\lambda_{ij}^*$, $i = 1, 2, \ldots , n-1,$ to the unit and
fixes other generators.
\end{lemma}

\begin{proof}
It is sufficient to prove that all defining relations turn to
the defining relations by such defined map. For the defining relations of $VP_{n-1}$
it is evident. If the relation of commutativity  (see relation (\ref{eq29}))
contains some generator of $V_{n-1}$
then by acting with $\varphi_n$ it turns to the trivial relation.
Let us consider the left hand side of relation (\ref{eq28}). We see  that it  contains
every
index two times.
Hence, if this part  includes some generator of $V_{n-1}$ (i.~e., one of the indices is equal
to $n$)
then some other generator contains the index $n$. Therefore, there
are
two generators of $V_{n-1}$ in the left hand side of the relation.
Since the right hand side contains all generators from the left hand side, then
by acting with $\varphi_n$ this relation turns to the trivial relation.
\end{proof}

\begin{lemma} \label{lemma5}
The following formulas are fulfilled in the group $VP_n$:\\

1)~~$ \lambda_{kl}^{\lambda_{ij}^{\varepsilon}} = \lambda_{kl},
~~~\mbox{max}\{i, j\} < \mbox{max}\{k, l\},~~~\varepsilon = \pm 1;
 $ \\

2)~~ $\lambda_{ik}^{\lambda_{ij}} = \lambda_{kj}^{\lambda_{ij}} \lambda_{ik}
\lambda_{kj}^{-1},~~~
\lambda_{ik}^{\lambda_{ij}^{-1}} = \lambda_{kj}^{-1} \lambda_{ik}
\lambda_{kj}^{\lambda_{ij}^{-1}},~~~i < j < k~~ \mbox{or}~~ j < i < k;
 $\\

3)~~$ \lambda_{ki}^{\lambda_{ij}} = \lambda_{kj} \lambda_{ki}
\lambda_{kj}^{-\lambda_{ij}},~~~
\lambda_{ki}^{\lambda_{ij}^{-1}} = \lambda_{kj}^{-\lambda_{ij}^{-1}}
\lambda_{ki} \lambda_{kj},~~~i < j < k~~ \mbox{or}~~ j < i < k;
 $\\

4)~~$ \lambda_{jk}^{\lambda_{ij}} = \lambda_{ik} \lambda_{jk}
\lambda_{kj} \lambda_{ik}^{-1} \lambda_{kj}^{-\lambda_{ij}},~~~
\lambda_{jk}^{\lambda_{ij}^{-1}} = \lambda_{jk}^{-\lambda_{ik}^{-1}}
\lambda_{ij}^{-1} \lambda_{jk} \lambda_{kj} \lambda_{ij},~~~i < j < k ~~ \mbox{or}~~ j < i <
k,
 $\\
where, as usual, different letters stand for different indices.
\end{lemma}

\begin{proof}
The formula 1) immediately follows from the first relation of
Theorem
\ref{theorem1}.

Consider relation (\ref{eq28}) from Theorem \ref{theorem1}:
$$
\lambda_{ki} \, (\lambda_{kj} \, \lambda_{ij}) = (\lambda_{ij} \, \lambda_{kj}) \, \lambda_{ki}.
$$
Note that the indices of generators are connected by one of the inequalities:
$$
a)~~ k < j < i,~~~b)~~ j < k < i,~~~c)~~ i < j < k,
$$
$$
d)~~ j < i < k,~~~e)~~ k < i < j,~~~f)~~i < k < j.
$$
If the indices are connected by inequality $a)$ or $b)$ then from (\ref{eq28}) we obtain
$$
\lambda_{ki}^{\lambda_{kj}} = \lambda_{ij}^{\lambda_{kj}} \,
\lambda_{ki} \,
\lambda_{ij}^{-1},
$$
and it is the first formula from
2).

If the indices in  relation (\ref{eq28}) are connected by inequality c) or d) we
obtain
$$
\lambda_{ki}^{\lambda_{ij}} = \lambda_{kj} \, \lambda_{ki} \,
\lambda_{kj}^{-\lambda_{ij}},
$$
and it is the first formula from 3).

If  indices in  relation (\ref{eq28}) are connected by inequality e)
or f) then
$$
\lambda_{ij}^{\lambda_{ki}} = \lambda_{kj} \, \lambda_{ij} \,
\lambda_{kj}^{-\lambda_{ki}}.
$$
Using the formula from 2), we obtain
$$
\lambda_{ij}^{\lambda_{ki}} = \lambda_{kj} \, \lambda_{ij} \, \lambda_{ji} \,
\lambda_{kj}^{-1} \,
\lambda_{ji}^{-\lambda_{ki}},
$$
and it is the first formula from 4).

The formulas of conjugations by elements $\lambda_{ij}^{-1}$ can be established similarly.

\end{proof}

Assume that the theorem is proven for the group $VP_{n-1}$. Hence, any element $w \in VP_{n-1}$
can be written in the form
$$
w \, = \, w_1 \, w_2 \ldots w_{n-2},~~~w_i \in V_i^*,
$$
where each word $w_i$ is a reduced word over the alphabet consisting if generators $\lambda_{ki}^{\pm 1}$, $1 \leq k \leq i-1,$ and reduced powers of generators
$\lambda_{ki}$, $1 \leq k \leq i-1,$ and their inverse.
Let us define  reduced powers of generators in the group $V_{n-1}^*$. We say that the element
$\lambda_{nk}(w) = \lambda_{nk}^w$ is {\it the reduced power of the generator}
$\lambda_{nk}$ if
$w$ is the empty word or a word written in the normal form and begin with
a
reduced power of some generator $\lambda_{lk}$ or its inverse.

The statement about decomposition in to the semi--direct product
$VP_n = V_n^* \leftthreetimes VP_{n-1}$ is
quite  evident. It remains to find generators of $V_n^*$ and prove its freedom.

\begin{lemma} \label{lemma6}
The group $V_{n-1}^*$ is a free group. It is generated by
$\lambda_{1n}, \lambda_{2n},$ $\ldots, $ $\lambda_{n-1,n}$ and all reduced powers of the
generators $\lambda_{n1}, \lambda_{n2}, \ldots, \lambda_{n,n-1}$.
\end{lemma}

\begin{proof}
The proof is similar to  that of Lemma \ref{lemma3}.
From Lemma \ref{lemma5} it follows that this set is the set of generators of
$V_{n-1}^*$.
Further, since the set of defining relations of $VP_n$ is equivalent to the set
of conjugating formulas  defining the action of $VP_{n-1}$ on $V_{n-1}^*$,
only trivial relations are fulfilled in $V_{n-1}^*$.
\end{proof}

Theorem \ref{theorem2}  follows from these results.

As a consequence of this theorem we obtain  the normal form of words in $VB_n$.

\begin{cor} \label{cor2}
Every element from $VB_n$ can be written uniquely in the form
$$
w = w_1 \, w_2 \ldots w_{n-1} \, \lambda,~~~\lambda \in
\Lambda_n,~~~w_i \in V_i^*,
$$
where $w_i$ is a reduced word in generators, reduced powers of generators and their inverse.
\end{cor}

The defined above homomorphism of the virtual braid group onto the welded braid group
agrees with the decomposition from Theorem
\ref{theorem2} and with the decomposition of $C_n \simeq WB_n$ described in the first section.

\begin{cor} \label{cor3}
The homomorphism $\varphi_{VW} : VB_n \longrightarrow WB_n$ agrees with the
decomposition of these groups,
i.~e., it maps the group $VP_n$ onto $Cb_n \simeq WP_n$ and the factors  $V_i^*$
 onto the factors $D_i$, $i = 1, 2, \ldots, n-1$.
\end{cor}

\section{The universal braid group}

Let us define  {\it the universal braid group} $UB_n$ as the group
with generators $\sigma_1, \sigma_2, \ldots, \sigma_{n-1},$
$c_1, c_2, \ldots, c_{n-1},$ defining relations
(\ref{eq1})--(\ref{eq2}), the relations:
$$
c_i \, c_j = c_j \, c_i,~~~ |i-j| \geq 2,
$$
and the mixed relations:
$$
c_i \, \sigma_j = \sigma_j \, c_i~~~ |i-j| \geq 2.
$$

Recall  (see \cite{BS}) that  Artin's group of the type $I$  is called the group
 $A_I$ with generators $a_i$, $i \in I$, and the defining relations
$$
a_i \, a_j \, a_i \ldots = a_j \, a_i \, a_j \ldots ,~~~ i, j \in I,
$$
where  words from the left and right hand sides consist of $m_{ij}$
alternating letters $a_i$ and $a_j$.

\begin{prop}\label{proposition1}
1) The group $UB_n$ has as a subgroup the braid group $B_n$.

2) There exist homomorphisms
$$
\varphi_{US} : UB_n \longrightarrow SG_n,~~~\varphi_{UV} : UB_n \longrightarrow
VB_n,~~~\varphi_{UB} : UB_n \longrightarrow B_n.
$$

3) The group $UB_n$ is  Artin's group.
\end{prop}

\begin{proof}
1) Evidently, there exists a homomorphism $B_n \longrightarrow UB_n$.
On the other hand, assuming
$\psi(\sigma_i) = \sigma_i$, $\psi(c_i) = e$, $i = 1, 2, \ldots, n-1,$ we obtain the
retraction
$\psi$ of $UB_n$ onto $B_n$. Therefore,
the subgroup $\langle \sigma_1, \sigma_2, \ldots, \sigma_{n-1} \rangle $ of
$UB_n$ is isomorphic to the braid group $B_n$.

2) Let us define the map $\varphi_{US}$  as follows
$$
\varphi_{US}(\sigma_i) = \sigma_i,~~~ \varphi_{US}(c_i) = \tau_i,~~~i = 1, 2, \ldots,
n-1.
$$
Comparing the defining relations of $UB_n$ and $SG_n$, we see that this map is a homomorphism.
Analogously, we can show that the map
$$
\sigma_i \longmapsto \sigma_i,~~~ c_i \longmapsto \rho_i,
$$
is extendable to the homomorphism $\varphi_{UV}$ and the map
$$
\sigma_i \longmapsto \sigma_i,~~~ c_i \longmapsto e,
$$
is extendable to the homomorphism $\varphi_{UB}$.

3) immediately follows from the defining relations
of $UB_n$ and the definition of Artin's group.
\end{proof}

It should be noted that none of the groups $SG_n$, $VB_n$, $WB_n$ (in the natural presentations)
is not Artin's group.

The following questions naturally arise in the context of the results obtained
above.

{\bf Problems.} 1) Solve the word and conjugacy problems in $UB_n$,
$n > 2$.\\
2) Is it possible to give some geometric interpretation for elements of $UB_n$
similar to the geometric interpretation for elements of the braid groups $B_n$, $SG_n$, $VB_n$, $UB_n$?










\begin{thebibliography}{99}

\bibitem{Vas}
V. A. Vassiliev, Complements of discriminants of smooth maps:
Topology and applications, Translations of Mathematical
Monographs, vol. 98, Amer. Math. Soc., Providence, RI, 1992.

\bibitem{Bir93}
J. S. Birman, New points of view in knot theory,
Bull. Am. Math. Soc., New Ser. 28, No.2 (1993), 253-287.


\bibitem{Ka}
L. H. Kauffman, Virtual knot theory,
Eur. J. Comb., 20, No.7 (1999), 663-690.

\bibitem{GPV} M. Goussarov, M. Polyak, O. Viro, Finite-type invariants of classical and
virtual knots,
Topology, 39, No.5 (2000), 1045-1068.


\bibitem{FRR}
R. Fenn, R. Rim\'{a}nyi, C. Rourke, The braid--permutation group,
Topology, 36, No.1 (1997), 123-135.

\bibitem{Baez}
J. C. Baez, Link invariants of finite type and perturbation
theory, Lett. Math. Phys., 26, No.1 (1992), 43-51.

\bibitem{Ver01}
V. V. Vershinin, On homology of virtual braids and Burau
representation, J. Knot Theory Ramifications, 10, No.5 (2001),
795-812.



\bibitem{Bir}
J. S. Birman, Braids, links and mapping class group,
Princeton--Tokyo: Univ. press, 1974.

\bibitem{Ge}
B. Gemein, Singular braids and Markov's theorem,
J. Knot Theory Ramifications, 6, No.4 (1997), 441-454.

\bibitem{Kam}
S. Kamada, Braid presentation of virtual knots and welded knots,
Preprint (math.GT/0008092).

\bibitem{Dash}
O. T. Dasbach, B. Gemein, The word problem for the singular braid
monoid, Preprint, 1999.

\bibitem{Corran}
R. Corran, A normal form for class of monoids including the
singular braid monoids, J. Algebra, 223, No.1 (2000), 256-282.

\bibitem{Ver}
V. V. Vershinin, On the singular braid monoid,
Preprint (math.GR/0309339).


\bibitem{Kr}
M. Guti\'{e}rrez, S. Krsti\'{c},
Normal forms for basis-conjugating automorphisms of a free group.
Int. J. Algebra Comput., 8, No. 6 (1998), 631-669.

\bibitem{Bar}
V. G. Bardakov,
The structure of the group of conjugating automorphisms,
Algebra i Logik, 42, No.5 (2003), 515-541.


\bibitem{BarP}
V. G. Bardakov,
The structure of the group of conjugating automorphisms and the linear representation
of the braid groups of some manifolds, Preprint (math.GR/0301247).


\bibitem{Mar}
A. A. Markoff, Foundations of the algebraic theory of braids, Trudy Mat. Inst. Steklova,
No. 16 (1945), 1--54.

\bibitem{Sav}
A. G. Savushkina,
On the group of conjugating automorphisms of a free group. (Russian, English)
Math. Notes, 60, No.1 (1996), 68-80; translation from Mat. Zametki 60, No.1 (1996), 92-108.

\bibitem{Mac}
J.~McCool, On basis--conjugating automorphisms of free groups, Can. J. Math., 38,
No. 6 (1986), 1525--1529.

\bibitem{FKR}
R. Fenn, E. Keyman, C. Rourke, The singular braid monoid embeds in a group,
J. Knot Theory Ramifications, 7, No.7 (1998), 881-892.


\bibitem{KMS}
W. Magnus, A. Karrass, D. Solitar, Combinatorial group theory,
Interscience Publishers, New York, 1996.




\bibitem{BS}
E. Brieskorn, K. Saito, Artin--Gruppen und Coxeter--Gruppen,
Invent. Math., 17 (1972), 245-271.

\end{thebibliography}
\end{document}